\documentclass[12pt]{article}
\usepackage{amssymb,amsmath}

\input epsf.tex

\title{Topological Borsuk problem}

\author {Yan Soibelman}

\begin{document}
\maketitle

\newtheorem{defn}{Definition}
\newtheorem{thm}{Theorem}
\newtheorem{lmm}{Lemma}
\newtheorem{rmk}{Remark}
\newtheorem{prp}{Proposition}
\newtheorem{conj}{Conjecture}
\newtheorem{exa}{Example}
\newtheorem{cor}{Corollary}
\newtheorem{que}{Question}
\newtheorem{ack}{Acknowledgements}
\newcommand{\K}{{\bf k}}
\newcommand{\C}{{\bf C}}
\newcommand{\R}{{\bf R}}
\newcommand{\N}{{\bf N}}
\newcommand{\Z}{{\bf Z}}
\newcommand{\Q}{{\bf Q}}
\newcommand{\G}{\Gamma}
\newcommand{\A}{A_{\infty}}
\newcommand{\ihom}{\underline{\Hom}}
\newcommand{\epi}{\twoheadrightarrow}
\newcommand{\mono}{\hookrightarrow}
\newcommand\ra{\rightarrow} 
\newcommand\uhom{{\underline{Hom}}}
\renewcommand\O{{\cal O}} 
\newcommand\nca{nc{\bf A}^{0|1}}

{\it To the memory of my wife Tanya with love}

\section{Introduction}

This paper is based on my  lecture at the US-Canada Mathcamp-2002
in Colorado Springs.
The text is very elementary (it was designed for advanced high school
students). At the same time the main problem about
the topological Borsuk number for ${\bf R}^n, n>2$ remains open
since 1977.

 I thank the organizers of the Mathcamp-2002 for inviting me to give
lectures.

\section{Borsuk conjecture and topological Borsuk number}

Famous Borsuk conjecture says that any compact $F$ in $\R^n$ can be partitioned
into $n+1$ closed subsets of smaller diameter.  The conjecture is true for $n=2$ and $n=3$ 
as well as  for all $F$ having smooth boundary. 
The latter result is topological. It is based on the following
Borsuk theorem: if the standard sphere $S^{n-1}\subset \R^n$ is represented
as a union $F_1\cup F_2\cup...\cup F_n$ of $n$ closed subsets then at least
one $F_i$ contains a pair of antipodal points $x$ and $-x$.

General Borsuk conjecture was disproved in [KK]. 
Let us call $b_{\R^n}(F)$ the minimal number of parts of smaller 
diameter necessary to
partition $F$. Then Kahn and Kalai constructed such $F$  that

$$b_{\R^n}(F)>(1.2)^{\sqrt{n}}.$$

The counterexample (as well as many of its simplifications)
 has combinatorial nature and deals with specific properties of the Euclidean
metric in $\R^n$.

Topological version of the Borsuk problem was formulated around 1977
(see [So]). We recall it below.

Let $(X,\rho_0)$ be a metric space ($\rho_0$ is the metric).
For any compact $F\subset X$ we denote by $b_{(X,\rho_0)}(F):=b_X(F)$
 its Borsuk
number, i.e. the minimal number of parts of smaller diameter necessary to partition $F$.
We denote by $B(X,\rho_0)$ the Borsuk number of the metric space $X$.
By definition, any compact in $X$ can be partitioned into $B(X,\rho_0)$ parts 
of smaller diameter, but there exist compacts which cannot be partitioned
into  $B(X,\rho_0)-1$ such parts.

Let $\Omega(\rho_0)$ be the set of metrics on $X$ 
which define the topology equivalent
to the one given  by $\rho_0$. We will call elements of this set
 $\rho_0$-equivalent metrics.

The number  $B(X,\rho)$ can change as 
$\rho$ varies inside of $\Omega(\rho_0)$.

\begin{exa} Let $\rho_0$ be the standard Euclidean metric
in $\R^2$ and $\rho$ be the Minkowski metric, i.e. $\rho((x_1,y_1),(x_2,y_2))=
|x_1-x_2|+|y_1-y_2|$ in coordinates. Clearly $\rho\in \Omega(\rho_0)$. 
Then the classical result, which goes back to 1950's, says that
$B(X,\rho)=4$. On the other hand $B(X,\rho_0)=3$.

\end{exa}

The following definition was  given in  [So].

\begin{defn} Topological Borsuk number of $(X,\rho_0)$ is defined as

$$B(X)=min_{\rho\in \Omega(\rho_0)}B(X,\rho).$$

\end{defn}

{\bf Topological Borsuk Problem}. {\it Estimate $B(X)$ for the Euclidean space $X={\bf R}^n$}.

In particular, is it true that $B(\R^n)\ge n+1$? More speculatively, is it true
that $B(\R^n)$ is bounded from below by $B(\R^n,\rho_0)$, where $\rho_0$ is the standard
Euclidean metric?

\section{2-dimensional case}

In order to prove that $B(X)\ge m$ it suffices to prove that for any metric
$\rho\in \Omega(\rho_0)$ there exists $c>0$ and a finite subset  $I\subset X$
consisting of $m$ elements such that $\rho(i,j)=c$  for all $i\ne j$.
If $X=\R^n$ and $\rho=\rho_0$ is the standard Euclidean metric then one
can take $c=1, m=n+1$ and $I$ be the set of vertices of the regular simplex.
It is not obvious that such $I$ exists for other $\rho_0$-equivalent metrics.
In the case $n=2$ the answer is positive due to the following theorem.

\begin{thm}([So]) Topological Borsuk number of the Euclidean $\R^2$ is equal to $3$.

\end{thm}

The proof presented below is basically the same as in [So].

{\it Proof}. Let $\rho$ be a metric on $X=\R^2$ which defines the topology equivalent to
the standard one. Let us consider the map $f:X^3\to \R^3$ such that

$$f(x_1,x_2,x_3)=(\rho_{12},\rho_{23},\rho_{13}),$$
where $\rho_{ij}=\rho(x_i,x_j), i\ne j$.

By our assumption the map $f$ is continuous. It suffices to prove that the image of $f$
intersects the line 
$l=\{(\rho_{12},\rho_{23},\rho_{13})|\rho_{12}=\rho_{23}=\rho_{13}\}$
 besides the obvious point $(0,0,0)$ which is the image of the diagonal
$x_1=x_2=x_3$.

Notice that $f$ is equivariant with respect to the natural actions of the
cyclic
group $\Z/3$ on $X^3$ and $\R^3$ (view $\Z/3$ as the subgroup 
of the symmetric group
$\Sigma_3$ acting by permutations of coordinates).  

It suffices to prove
 that there is no continuous
 $\Z/3$-equivariant map between $X^3\setminus f^{-1}(0)$
and $\R^3\setminus l$.

Since $X^3=\R^6$ and $f^{-1}(0)=\{(x_1,x_2,x_3)|x_1=x_2=x_3\}$ 
the set $X^3\setminus f^{-1}(0)$ is equivariantly homotopic
to the $3$-dimensional sphere $S^3$ (which can be considered as
 a $\Z/3$-invariant subset 
 of the plane $x_1+x_2+x_3=0$).
Similarly, $\R^3\setminus l$ is equivariantly homotopic to the 
$1$-dimensional sphere $S^1$ 
considered
as a $\Z/3$-invariant 
subset of the plane $\rho_{12}+\rho_{23}+\rho_{13}=0$. Notice
that the natural action of $\Z/3$ is free on both spheres.
It is enough to prove that there is no $\Z/3$-equivariant map between $S^3$ and $S^1$.
It was done in [So] by using the notion of the category (genus) of a topological space.
Here we will use its modern version called the $G$-index of a topological space
 where $G$ is a group acting on the space
(see [M]). In our case $G=\Z/3$.

\begin{defn} Let $G$ be a non-trivial  finite group. An $E_nG$-space is a $G$-space
$Y$ such that

a) $G$ acts freely on $Y$;

b) $Y$ is $n$-dimensional;

c) $Y$ is $(n-1)$-connected.

\end{defn}

If $X$ is a $G$-space then $ind_G(X)$ (the $G$-index of $X$) is the minimal $n$ such
that there exists a $G$-equivariant map $X\to E_nG$ (the index can be infinite).It is easy to see that $ind_G(E_nG)=n$.
One can prove (see [M], 6.2.5) that there is no $G$-equivariant continuous map
$E_nG\to E_{n-1}G$ (Borsuk-Ulam type theorem). More generally, there is no
$G$-equivariant continuous map $f:X\to Y$ between $G$-spaces $X$ and $Y$ such
that $ind_G(X)>ind_G(Y)$.

 On the other hand, if $p$ is a prime
number then any odd-dimensional sphere is a $\Z/p$-space. Indeed, the group
$\Z/p$ acts on $S^{2n-1}=\{(z_1,...,z_n)\in \C^n|\sum_i|z_i|^2=1\}$ via
$(z_1,...,z_n)\mapsto (z_1exp(2\pi i/p),...,z_nexp(2\pi i/p))$.
 Taking $p=3$ we
finish the proof of the theorem. $\blacksquare$

\section{Conclusion}

The proof above uses only elementary algebraic topology
and very little information about the metric. The proof does not work when we
replace $\R^2$ by $\R^n, n>2$. In that case we have a 
$\Sigma_{n+1}$-equivariant continuous
map $(\R^n)^{n+1}\to R^{n(n-1)/2}$.
 Although arising topological spaces are 
still spheres, the natural actions of various subgroups of the symmetric group
$\Sigma_{n+1}$ are not free on the target.
Perhaps one needs new ideas in order to estimate the topological Borsuk number of $\R^n$.

\vspace{10mm}

{\bf References}

\vspace{5mm}

[KK] J. Kahn, G. Kalai, A counterexample to Borsuk's conjecture. Bull. AMS, 29(1993), 60-62.

[M] J. Matousek, Topological Methods in Combinatorics and Geometry (book, to appear).

[So] Y. Soibelman, Calculation of a topological invariant
connected with the partitioning of figures into parts of smaller diameter. (Russian) Mat.
Zametki 27 (1980), no. 4, 647--649, 671. English Translation in Math. Notes 27 (1980), 3-4,
317-318.

\vspace{10mm}

Address: Department of Mathematics, KSU, Manhattan, KS 66506, USA

{soibel@math.ksu.edu}

\end{document}